\documentclass[review,12pt]{elsarticle}

\usepackage{graphicx}
\usepackage{amsfonts}
\usepackage{amsmath}
\usepackage{amssymb}
\usepackage{multirow}
\usepackage{color}
\usepackage[version=4]{mhchem}

\usepackage{siunitx}
\graphicspath{{plots/}}

\journal{Journal of Computational Physics}

\begin{document}

\begin{frontmatter}

\title{A finite element method for stochastic diffusion equations using fluctuating hydrodynamics}

\author{P. Martínez-Lera\corref{firstcorr}}
\cortext[firstcorr]{p.martinez@unizar.es}

\author{M. De Corato\corref{secondcorr}}
\cortext[secondcorr]{mdecorato@unizar.es}

\address{Aragon Institute of Engineering Research (I3A), University of Zaragoza, \\50018 Zaragoza, Spain}

\begin{abstract}
We present a finite element approach for diffusion problems with thermal fluctuations based on a fluctuating hydrodynamics model. The governing equations are stochastic partial differential equations with a fluctuating forcing term. We propose a discrete formulation of the fluctuating forcing term that has the correct covariance matrix up to a standard discretization error. Furthermore, we derive a linear mapping to transform the finite element solution into an equivalent discrete solution that is free of the artificial correlations introduced by the spatial discretization. 
The method is validated by applying it to two diffusion problems: a second-order diffusion equation and a fourth-order diffusion equation. The theoretical (continuum) solution to the first case presents spatially decorrelated fluctuations, while the second case presents fluctuations correlated over a finite length. In both cases, the numerical solution presents a structure factor that approximates well the continuum one.
\end{abstract}


\end{frontmatter}

\section{Introduction}
With microfluidics technologies becoming increasingly widespread and multiple applications already reaching the market level \cite{volpatti2014commercialization}, the great challenges of technology now shift towards the nanoscale, with potential for multiple breakthroughs and ground-breaking innovations. In analogy with microfluidics, the rapidly growing field of research that studies the flow and transport phenomena in fluidic environments of nanoscopic dimensions has been termed “nanofluidics” \cite{bocquet2020nanofluidics}. The blooming of nanofluidics was enabled in the last decade by the improvements and development of new experimental and computational techniques and has shown multiple phenomena that have no counterparts at larger scales \cite{faucher2019critical,laine2020nanotribology}.

Fluid flow, solute transport and chemical reactions look very different at the nanoscale than they do at the macroscopic level. At the small scales, the thermal fluctuations of solvent and solute molecules are relevant and cannot be ignored. Even in situations where one can average over many molecules, quantities like mass, temperature and momentum fluctuate around their mean value \cite{de2006hydrodynamic}. These fluctuations are very well understood around equilibrium \cite{de2006hydrodynamic} but can lead to unexpected and highly nontrivial phenomena when they occur in systems driven out of equilibrium \cite{donev2011diffusive,peraud2017fluctuation,vutukuri2020active}.

As most of the applications in energy harvesting, membrane technology and biomedicine involve systems driven out from equilibrium, understanding the effects of thermal fluctuations in these instances is of utmost importance. Thermal fluctuations are also particularly large inside eukaryotic cells, where some molecules are present in very small numbers and need to be recruited in specific locations to achieve their biological functions \cite{brown2013linking}.
Deterministic models fail to predict many of the important consequences of thermal fluctuations, which calls for the development of efficient numerical tools that can include them within an efficient framework.

To do so, one possibility is to modify the deterministic partial differential equations (PDEs) that are typically used to model phenomena at larger scales to include stochasticity.  
One such approach to incorporate the effect of the thermal fluctuations in the transport equations was proposed by Landau and Lifschitz \cite{landau2013statistical}, and is based on adding a stochastic forcing term that satisfies the fluctuation-dissipation balance. The transport equations that result from adding a fluctuating forcing term are stochastic partial differential equations (SPDEs).
This approach to include thermal fluctuations into deterministic PDEs has been called fluctuating hydrodynamics \cite{de2006hydrodynamic} and it has been used to study stochastic drift-diffusion equations
(e.g. 
\cite{donev2011diffusive,leonard2013stochastic,chaudhri2014modeling,peraud2016low,donev2019fluctuating}),
stochastic reaction-diffusion equations
(e.g. 
\cite{atzberger2010spatially,bhattacharjee2015fluctuating,kim2017stochastic}),
the Brownian motion of colloidal particles 
(e.g. 
\cite{usabiaga2013inertial,keaveny2014fluctuating,delmotte2015simulating,de2016finite,sprinkle2019brownian,westwood2022generalised}),  {fluctuations of viscoelastic fluids \cite{hutter2018fluctuating,hutter2020fluctuating},}
stochastic thin-film equations
(e.g. \cite{tsekov1993effect,sprittles2023rogue,liu2023thermal}),
diffusion within membranes and liquid-liquid interfaces
(e.g. \cite{wang2013dynamic,rower2022surface}), and bubble nucleation and growth
(e.g.
\cite{gallo2018thermally,gallo2023nanoscale}).
The advantage of a fluctuating hydrodynamics approach compared to discrete particle methods is apparent when one needs to study mesoscopic systems that involve a vast number of particles and long timescales compared to the characteristic diffusion time of individual particles.

In this work, we propose a finite element framework to solve diffusion equations using fluctuating hydrodynamics. The governing equations are SPDEs of the kind
\begin{equation}\label{eq:PDEintro}
	\frac{\partial {u}}{\partial t} - \nabla \cdot \left(D \nabla F(u)\right) = \nabla \cdot \left( \sqrt{2 D u}\boldsymbol{\zeta} \right)\,.
\end{equation}
where $u$ is the number density of a diffusing species, $D$ is its diffusion coefficient and
$\boldsymbol{\zeta}$ represents white noise in space and time, with stochastic properties given by
\begin{equation}\label{eq:meanwhitenoise}
\left<\zeta_i(\mathbf{x},t)\right> = 0 \,
\end{equation}
and
\begin{equation}\label{eq:varwhitenoise}
\left<\zeta_i(\mathbf{x},t) \zeta_j(\mathbf{x}',t')\right> = \delta_{ij} \delta(\mathbf{x}-\mathbf{x}')\delta(t-t') \,.
\end{equation}
In the case of $F(u)=u$, Equation \eqref{eq:PDEintro} can be formally derived from the motion of independent point-sized Brownian walkers using Ito calculus \cite{dean1996langevin,donev2014dynamic}, and it is nothing more than a rewriting of the equations of motion of a collection of independent Brownian walkers in terms of their density, $u$.  {Stochastic equations of the form of Eq. \eqref{eq:PDEintro} also arise in the context of Dynamic Density Functional Theory \cite{donev2014dynamic,archer2004dynamical}. While an ensemble-average of the stochastic equations describing the trajectories of the Brownian walkers would lead to a deterministic PDE describing the time evolution of the particle distribution, an SPDE like Equation \eqref{eq:PDEintro} describes the evolution of a coarse-grained density distribution instead. Detailed discussions on the derivation of deterministic and stochastic equations in this context and their different interpretations can be found in \cite{archer2004dynamical,donev2014dynamic,chavanis2015generalized,chavanis2019generalized}.}

In  Equation \eqref{eq:PDEintro} we have abused notation to highlight its formal similarity to deterministic PDEs, but, in fact, neither the solution $u$ nor the white-noise source term $\boldsymbol{\zeta}$
can be interpreted as differentiable point-wise functions. Still, as mentioned above, numerical solutions to equations based on fluctuating hydrodynamics such as Equation \eqref{eq:PDEintro} are widely used to capture the physical effects from the effect of thermal fluctuations.
Furthermore, recent advances in the theoretical analysis front seem to confirm that the truncation of the fluctuations for small wavelengths, including the natural truncation resulting from standard spatial discretization techniques, leads to well-posed fluctuating hydrodynamics equations \cite{cornalba2023dean}.

When solving such equations, one is interested in the expected value of the solution, but also in its second-order statistical moments. For instance, the so-called structure factor, related to the spatial Fourier transform of the autocorrelation function of the solution, can be measured experimentally \cite{de2006hydrodynamic}.  
However, both the choice of spatial discretization and temporal integration schemes influence the fluctuation-dissipation balance and can introduce artificial spatial correlations in the solution, thus leading to non-physical second-order statistical moments. \cite{delong2013temporal,donev2010accuracy} 

For this reason, the effect of time integrators on the solution of fluctuating-hydrodynamics equations has been investigated in depth \cite{delong2013temporal,donev2010accuracy}. 
As for the spatial discretization, finite volume methods 
(e.g. \cite{donev2010accuracy,balboa2012staggered,russo2021finite, kim2017stochastic}) can be designed so as not to introduce additional artificial correlations, since each degree of freedom typically corresponds to a single finite volume, for which the solution can be integrated independently of the rest of finite volumes. This is not the case, however, for finite element methods in general, which yield solutions with artificial correlations that depend on the choice of shape functions \cite{de2015finite}. This makes finite element solutions difficult to interpret physically and nonlinear terms in the equations difficult to handle mathematically. 
To overcome this limitation for general finite element methods, it has been suggested that specialized finite element discretization schemes should be used, which are designed based on physical considerations instead of just numerical analysis \cite{de2015finite}. The rationale behind this strategy is based on the idea that discretization and physical coarse-graining are deeply connected. 

In the present work, we start from a different perspective and propose an approach that yields a finite element solution to Equation \eqref{eq:PDEintro} in which both the expected value of the solution and its second-order statistical moments are well approximated. In particular, we propose to first obtain a numerical solution using any suitable spatial discretization that is selected based on numerical analysis only and then apply a linear transformation that removes from the numerical solution the artificial correlations introduced by the chosen spatial discretization. 
With the proposed approach, the choice of spatial discretization techniques is not restricted anymore to those that produce any particular spatial correlations, and we can use the finite element method and its powerful mathematical framework to deal with complex geometries and boundary conditions without limitations. The resulting solution has a physical interpretation that is detached from the numerical discretization, just as numerical solutions to deterministic PDEs do.

We work in a finite element framework with a standard Galerkin discretization. Nevertheless, our conclusions are generalizable to other numerical approaches that involve a spatial discretization. We apply our method to a second-order diffusion equation and to a fourth-order diffusion equation. The two test cases are complementary, since one presents only short-ranged fluctuations, while the other has a solution that is spatially correlated over a finite length. Both cases were already used in Ref. \cite{de2015finite} to test a Petrov-Galerkin finite element implementation. 

The outline of this paper is as follows: In Section \ref{sec:FE} we describe the finite element discretization for a second-order fluctuating diffusion equation that has a solution with short-ranged correlations. We also propose an implementation for the fluctuating forcing term that has the correct covariance up to a normal discretization error. In Section \ref{sec:results_diff2}, we present finite element results for a one-dimensional boundary value problem of the second-order equation, and we discuss the presence of artificial spatial correlations in the solution due to the effects of the spatial discretization and temporal integration.
In Section \ref{sec:decorrelation}, we propose a linear mapping to remove from the numerical solution the artificial correlations introduced by the spatial discretization, and we apply it to the results of the same second-order diffusion equation.
In Section \ref{sec:results_diff4}, we apply the proposed method to a one-dimensional fourth-order diffusion equation, to test its performance in the case of a solution that is spatially correlated over a finite length. Finally, we summarize our conclusions in Section \ref{sec:conclusions}.

%%%

\section{Numerical approximation}\label{sec:FE}

In this section, we describe the numerical schemes used to solve Equation \eqref{eq:PDEintro} numerically. We focus here on the case 
$F(u)=u$,
for which numerical results will be presented in Sections \ref{sec:results_diff2} and \ref{sec:decorrelation}. Nevertheless, most of the discretization details presented in this section are also relevant for the case  $F(u)=u +\left(\ell_0/2\pi\right)^2 \nabla^2 u\,,$ %$$F(u)=u +\left(\frac{\ell_0}{2\pi}\right)^2 \nabla^2 u\,,$$
which will be discussed in Section \ref{sec:results_diff4}.
 
\subsection{Weak form and finite element discretization}\label{sec:discretiz}
We consider the second-order equation
\begin{equation}\label{eq:diff2}
	\frac{\partial {u}}{\partial t} - D \nabla^2 u = \nabla \cdot \left( \sqrt{2 D u}\boldsymbol{\zeta} \right)\,
\end{equation}
with constant diffusion coefficient $D$. This corresponds to Dean's equation \cite{dean1996langevin}, which describes the fluctuations of the density in a system of particles undergoing independent Brownian motions. The concentration $u$ represents therefore a coarse-grained number of particles per unit volume. 

A weak formulation of Equation \eqref{eq:diff2} can be obtained by multiplying it by the conjugate of a test function $v$ 
and integrating over the domain $\Omega$. After integration by parts, one obtains
\begin{multline}\label{eq:diff2_weak}
	\int_{\Omega}{v^*\frac{\partial {u}}{\partial t}}\,d^3\mathbf{x} + \int_{\Omega}{D \nabla v^* \cdot \nabla u}\,d^3\mathbf{x} = - \int_{\Omega}{\nabla v^* \cdot \left( \sqrt{2 D u}\boldsymbol{\zeta} \right)}\,d^3\mathbf{x} +
\\+ \int_{\partial \Omega}{v^* \left( \sqrt{2 D u}\boldsymbol{\zeta} \right)\cdot \mathbf{n}}\,d^2\mathbf{x}+
 \int_{\partial \Omega}{D v^* \nabla u \cdot \mathbf{n}}\,d^2\mathbf{x}
 \,,
\end{multline}
where $\partial \Omega$ is the boundary of the domain and $\mathbf{n}$ its normal unitary vector.

The domain $\Omega$ is discretized by dividing it into non-overlapping elements $\Omega_e$, such that
$\Omega=\cup_e \Omega_e$ with
$\Omega_e \cap \Omega_{e'} = \varnothing$ for $e \neq e'$. The solution space, as well as the test function space, are discretized using a given basis of shape functions. Here we use a standard Galerkin discretization with Lagrange polynomials as shape functions $\phi_i=\phi_i(\mathbf{x})$ for both the test function $v$ and the solution $u$, such that
\begin{equation}\label{eq:discreteu}
u(\mathbf{x}, t) \approx \sum_{i=1}^{N_{\text{dof}}}{\widetilde{u}_i(t) \phi_i(\mathbf{x})} \,,
\end{equation}
and
\begin{equation}\label{eq:discretev}
v(\mathbf{x}, t) \approx \sum_{i=1}^{N_{\text{dof}}}{\widetilde{v}_i(t) \phi_i(\mathbf{x})} \,.
\end{equation}
where $N_{\text{dof}}$ is the total number of degrees of freedom. 

The goal of the
finite element method is to find the values $\widetilde{u}_i(t)$ that satisfy the discretized weak equation \eqref{eq:diff2_weak} for any values $\widetilde{v}_i(t)$. Using Equations \eqref{eq:diff2_weak}, \eqref{eq:discreteu} and \eqref{eq:discretev}, we can express this problem as a linear system of equations in matrix form
\begin{equation}\label{eq:matricial2}
\mathbf{M} \frac{\partial \widetilde{\mathbf{u}}}{\partial t} + \mathbf{D}  \widetilde{\mathbf{u}} = \mathbf{f}+\mathbf{f}_{BC} \,,
\end{equation}
where $\mathbf{M}$ is the so-called mass matrix
\begin{equation}\label{eq:massmatrix}
M_{ij} = \int_{\Omega}{\phi_i \phi_j}d^3\mathbf{x} \,
\end{equation}
$\mathbf{D}$ is the diffusion matrix
\begin{equation}
D_{ij} = \int_{\Omega}{D \nabla \phi_i \cdot \nabla \phi_j}d^3\mathbf{x} \,,
\end{equation}
$\widetilde{\mathbf{u}}$ is an array that contains the coefficients $\widetilde{u}_i$ for each degree of freedom, $\mathbf{f}_{BC}$ is an array that incorporates the effect of the boundary conditions in the system and $\mathbf{f}$ is the random excitation term defined such that 
\begin{equation}\label{eq:fi}
f_i(t) = -  \int_{\Omega}{\nabla \phi_i \cdot \sqrt{2 D u} \mathbf{\boldsymbol{\zeta}}(t)}d^3\mathbf{x} \,.
\end{equation}
The covariance matrix of the components $f_i(t)$ can be expressed as 
\begin{equation}\label{eq:ficov}
\left<f_i(t) f_j(t')\right> = {2 D \delta(t-t')} \int_{\Omega}{ u(\mathbf{x}) \nabla \phi_i(\mathbf{x}) \cdot  \nabla \phi_j(\mathbf{x})}d^3\mathbf{x}\,,
\end{equation}
where $\left<.\right>$ indicates the expected value.
The implementation of $\mathbf{f}$, which represents the effect of the thermal fluctuations, is discussed in Section \ref{sec:noiseterm}, and the time integration schemes are presented in Section \ref{sec:timeintegrators}.

%%%

%%%
\subsection{Finite element formulation for the thermal fluctuations}
\label{sec:noiseterm}

To find a discrete approximation of $\mathbf{f}$, we seek to postulate a formulation that satisfies Equation \eqref{eq:ficov}. Here we propose an approach inspired by previous work in the context of the Stokes equations with fluctuating hydrodynamics \cite{de2016finite}. To this end, we make use of the same numerical integration schemes that are used in the finite element implementation, and that allow the integration of a given function $g(\mathbf{x})$ as
\begin{equation}
\int_{\Omega}g(\mathbf{x})d^3\mathbf{x}\approx \sum_{k=1}^{N_g}{w_k g(\mathbf{x}_k)}\,,
\end{equation}
where $N_g$ is the number of integration points, $g(\mathbf{x}_k)$ is the value of the function evaluated at the integration point and $w_k$ is the weight corresponding to integration point $k$. If an isoparametric formulation is used, $w_k$ would include both the weight corresponding to the integration rule in the reference element as well as the Jacobian for the mapping from the reference element to the physical element.
We assume that $w_k\ge0$. With this, we can postulate the following formulation for array $\mathbf{f}$ before the time discretization, as defined by Equation \eqref{eq:fi}:
\begin{equation}\label{eq:fi_integr}
    f_i=-\sqrt{2 D} \sum_{k=1}^{N_g}{\sqrt{w_k u(\mathbf{x}_k,t)} \boldsymbol{\zeta}_{k}(t) \cdot \nabla \phi_i(\mathbf{x}_k)}\,,
\end{equation}
where the components of $\boldsymbol{\zeta}_k$ are stochastic Gaussian process with $\left<\zeta_{k,m}\right>=0$ and  $\left<\zeta_{k,m}(t) \zeta_{l,m}(t')\right>=\delta(t-t') \delta_{kl}$.
The formulation given by Equation \eqref{eq:fi_integr} leads to the following covariance matrix
\begin{equation}\label{eq:fflu0}
\begin{split}
\left<f_i\right.&\left.(t)  f_j(t')\right>  = \\ \,&={2 D} \sum_{m=1}^{3}{\sum_{k=1}^{N_g}{\sum_{l=1}^{N_g}{\sqrt{w_k w_l u(\mathbf{x}_k,t) u(\mathbf{x}_l,t')} \frac{\partial \phi_i(\mathbf{x}_k)}{\partial x_{m}} \frac{\partial \phi_j(\mathbf{x}_l)}{\partial x_{m}} \left<\zeta_{k,m}(t) \mathbf{\zeta}_{l,m}(t') 
\right>}}}\\ \,
 &= {2 D \delta(t-t')} \sum_{k=1}^{N_g}{{u(\mathbf{x}_k,t)} \nabla \phi_i(\mathbf{x}_k) \cdot \nabla \phi_j(\mathbf{x}_k) {w_k}}
\,.
\end{split}
\end{equation}
This covariance matrix is identical to that of Equation \eqref{eq:ficov} up to an error introduced by the numerical integration scheme, thus proving that the postulated formulation for $\mathbf{f}$ in Equation \eqref{eq:fi_integr} has the correct second-order statistical moments. To preserve the fluctuation-dissipation balance, the integration rule used to compute the diffusion matrix $\mathbf{D}$ and the forcing term $\mathbf{f}$ should be the same \cite{de2016finite}.
In the particular case that the value of the concentration $u$ is large enough to linearize the fluctuating forcing term, the co-variance of the stochastic noise term  for degrees of freedom $i$ and $j$ becomes
\begin{equation}\label{eq:cov_b}
	\left< f_i(t) f_j(t') \right> = 2 u_0\, \delta(t-t') {D}_{ij} \,
\end{equation}
where $D_{ij}$ is the $(i,j)$ component of the diffusion matrix $\mathbf{D}$ and $u_0$ represents the average concentration. In \ref{sec:appendixfluc}, an alternative to Equation \eqref{eq:fi_integr} to define the linearized $\mathbf{f}$ is described, based on a decomposition of matrix $\mathbf{D}$. This alternative implementation of the linearized forcing term has computational advantages and is equivalent to Equation \eqref{eq:fi_integr} if the average concentration $u_0$ in the domain is large enough.

%%%

\subsection{Time integration and fluctuation-dissipation balance}
\label{sec:timeintegrators}
The effect of the temporal integrators for fluctuating-hydrodynamics equations has been explored extensively \cite{donev2010accuracy,delong2013temporal}. In this work, we use one-stage time integration schemes. For Equation \eqref{eq:matricial2}, the result for each time iteration can be expressed as
\begin{equation}
\widetilde{\mathbf{u}}^{n+1} = \left[\mathbf{M} + (1-\alpha) \Delta t\mathbf{D}\right]^{-1} \left[ \left(\mathbf{M} - \alpha \Delta t \mathbf{D}\right)\, \widetilde{\mathbf{u}}^{n} + \sqrt{\Delta t}\, \mathbf{f}^n + \Delta t \mathbf{f}_{BC}^n \right] \,,
\label{eq:iterator}
\end{equation}
where the superindex $(.)^n$ indicates that the term is evaluated at the time step $n$, and $\alpha$ is a constant that indicates whether the time integration is explicit ($\alpha=1$), implicit ($\alpha=0$) or semi-implicit ($\alpha=1/2$). Based on Equation \eqref{eq:fi_integr} in time, the stochastic forcing term in Equation  \eqref{eq:iterator} can be expressed as
\begin{equation}\label{eq:fi_integr_timeint}
    f_i^n=-\sqrt{2 D} \sum_{k=1}^{N_g}{\sqrt{w_k u(\mathbf{x}_k,t_n)} \boldsymbol{z}^n_{k} \cdot \nabla \phi_i(\mathbf{x}_k)}\,,
\end{equation}
where the components of $\boldsymbol{z}_k$ are stochastic Gaussian process with $\left<z_{k,m}\right>=0$ and  $\left<z_{k,m}(t) z_{l,m}(t')\right>= \delta_{kl}$.
We note that the solution, $u(\mathbf{x}_k,t_n)$, at the integration point, $\mathbf{x}_k$, has to be evaluated at time $t_n$ because the stochastic term has been derived using an Ito interpretation \cite{dean1996langevin,donev2014dynamic,kloeden1992stochastic}.

If the fluctuating term can be linearized, the alternative expression provided by Equation \eqref{eq:I_linear} can be used.
The elements of the covariance matrix of the discrete forcing term given by Equation \eqref{eq:fi_integr_timeint} can be expressed as 
\begin{equation}\label{eq:fflu0_timeint}
\left<f_i^n f_j^m\right>  = 
  {2 D \delta_{nm}} \sum_{k=1}^{N_g}{{u(\mathbf{x}_k,t)} \nabla \phi_i(\mathbf{x}_k) \cdot \nabla \phi_j(\mathbf{x}_k) {w_k}}
\,.
\end{equation}

In Equation \eqref{eq:iterator}, the value $\alpha=1/2$ corresponds to the implicit midpoint method (Crank-Nicholson scheme). This scheme is not only stable regardless of the value of the time step, but it also produces a solution with static spatial correlations that are independent of the time step as well \cite{donev2010accuracy,delong2013temporal}. 
Indeed, if we consider for now an infinite domain, so that we can neglect the effect of the boundary conditions, and we take the covariance of Equation \eqref{eq:iterator}, we obtain
\begin{equation}
\mathbf{D}\, \mathbf{C}\, \mathbf{M}^T +  \mathbf{M}\, \mathbf{C}\, \mathbf{D}^T 
+ \Delta t (1-2\alpha) \mathbf{D}\, \mathbf{C}\, \mathbf{D}^T
 = \mathbf{C}_{ff} \,,
\label{eq:covlhs}
\end{equation}
where $\mathbf{C}_{ff}$ is the covariance matrix of the forcing term $\mathbf{f}$, given by Equation \eqref{eq:fflu0_timeint}, and $\mathbf{C}$ is the covariance matrix of the discrete solution 
\begin{equation}
\mathbf{C}=\left<\left(\mathbf{\widetilde{u}}-\left<\mathbf{\widetilde{u}}\right>\right)\left(\mathbf{\widetilde{u}}^T-\left<\mathbf{\widetilde{u}}^T\right>\right)\right>\,.
\end{equation}
Equation \eqref{eq:covlhs} expresses the fluctuation-dissipation balance at the discrete level. The term that depends on the time step vanishes if $\alpha=1/2$.
By contrast, any other value of $\alpha$, including $\alpha=1$ (fully explicit) and $\alpha=0$ (fully implicit), will not make this term vanish and, as a result, will lead to a solution with spatial correlations that depend on the ratio   
\begin{equation}
\beta=\frac{D \Delta t}{\Delta x^2}\,,
\end{equation}
where $\Delta t$ and $\Delta x$ are the time step and the reference cell size, respectively. 

The choice of time integration scheme influences, therefore, the fluctuation-dissipation balance, and, therefore, the spatial correlations of the solution. In this work, we use the Crank-Nicholson scheme ($\alpha=1/2$) unless otherwise stated. With this scheme, Equation \eqref{eq:covlhs} becomes
\begin{equation}
\mathbf{D}\, \mathbf{C}\, \mathbf{M}^T +  \mathbf{M}\, \mathbf{C}\, \mathbf{D}^T 
 = \mathbf{C}_{ff} \,,
\label{eq:covlhs2}
\end{equation}
which shows that the fluctuation-dissipation balance is independent of the time step. In the linearized case, Equation \eqref{eq:covlhs2} becomes
\begin{equation}
\mathbf{D}\, \mathbf{C}\, \mathbf{M}^T +  \mathbf{M}\, \mathbf{C}\, \mathbf{D}^T 
 = 2 u_0 \mathbf{D} \,.
\label{eq:covlhs3}
\end{equation}
It is worth noting that in general, the fluctuation-dissipation balance is modified by the presence of boundary conditions, which we are not considering in this simplified analysis.

We see from Equation \eqref{eq:covlhs3} that, even with a choice of a temporal integrator that does not introduce additional correlations in the solution, the covariance matrix of the solution will depend on the mass matrix $\mathbf{M}$. As a result, the discrete solution may display spatial correlations that are not physical but induced by the spatial discretization.  {Indeed, if the diffusion and mass matrices are symmetric, the covariance matrix
\begin{equation}
\mathbf{C}=u_0 \mathbf{M}^{-1}
\end{equation}
satisfies Equation \eqref{eq:covlhs3}. If $\mathbf{M}$ is diagonal, $\mathbf{C}$ will also be diagonal, and the solution for each degree of freedom will be uncorrelated with the solution for the rest of degrees of freedom; conversely, if $\mathbf{M}$ is not diagonal, as is in general the case with finite-element methods because of the way the shape functions are defined in the elements, the solution for the different degrees of freedom will be correlated.}
It is worth noting that these artificial correlations do not represent a numerical error, but are instead a natural consequence of the choice of shape functions. Still, even if these artificial correlations do not represent a numerical error, they make the solution difficult to interpret and to compare to experimental measurements, and they can make non-linear terms of fluctuating-hydrodynamics equations difficult to treat.
In Section \ref{sec:results_diff2}, we show finite element results for a (bounded) 1D problem with periodic boundary conditions that present discretization-induced spatial correlations in the solution. Furthermore, in Section \ref{sec:decorrelation}, we propose a change of basis to transform the discrete solution into an equivalent discrete solution that is free of these artificial correlations.

\subsection{Coarse-graining and discretization}
\label{sec:coarsegraining}

Stochastic diffusion equations such as the ones that we are considering in this work can be seen as coarse-grained descriptions of the underlying system of particles. In fact, the macroscopic diffusion equations can be obtained from a description of the microscopic system through a coarse-graining procedure \cite{de2011coarse,espanol2015coupling}. 
In the case of diffusion equations such as Equation \eqref{eq:PDEintro}, the concentration $u$ is a statistical representation of the distribution of the particles moving in the domain. In particular, $u$ satisfies
\begin{equation}
    \int_{\Omega}{\gamma_i(\mathbf{x}) u(\mathbf{x},t)}d^3\mathbf{x}=\sum_{p=1}^{N_p}{\gamma_i(\mathbf{x}_p(t))}\,,
\end{equation}
where $N_p$ is the number of particles of the system, $\mathbf{x}_p$ their positions and $\gamma_i$ is a suitable coarse-graining function associated with the discrete coarse-graining volume $i$.

Using the approximation given by Equation \eqref{eq:discreteu}, we obtain
\begin{equation}\label{eq:coarse3}
    \int_{\Omega}{\gamma_i(\mathbf{x}) \sum_j \phi_j(\mathbf{x}) \widetilde{u}_j(t)} d^3\mathbf{x}=\sum_{p=1}^{N_p}{\gamma_i(\mathbf{x}_p(t))} =\widetilde{n}_i\,.
\end{equation}
In this work, we want to stress that, although the spatial correlations of the coarse-grained variable $\widetilde{\mathbf{n}}$ depend on the functions $\gamma_i$,  the spatial correlations of the discrete concentration $\widetilde{\mathbf{u}}$ depend on the choice of shape functions $\phi_i$ used to approximate $u$ as expressed by Equation \eqref{eq:discreteu}, and not on the coarse-graining functions $\gamma_i$. We will discuss this further in Section \ref{sec:decorrelationMatrix}, where we propose an approach to find a discrete approximation of the continuum concentration $u$ that is free of discretization-induced spatial correlations. 

%%%

\section{Finite element results for a 1D diffusion problem}
\label{sec:results_diff2}
In this section, we present numerical results for a one-dimensional version of Equation \eqref{eq:diff2}, using the discretization schemes described in Section \ref{sec:FE}, and illustrate the effects of the spatiotemporal discretization on the spatial correlations of the discrete solution.

\subsection{One-dimensional boundary value problem}
We consider a one-dimensional diffusion problem with thermal fluctuations, governed by equation
\begin{equation}\label{eq:diffusion}
	\frac{\partial u({x},t)}{\partial t} -  D \frac{\partial^2 u({x},t)}{\partial x^2}   = \frac{\partial}{\partial x} \left( \sqrt{2 D u} \mathbf{{\zeta}}({x},t)  \right)\,,
\end{equation}
where $D$ is the diffusion coefficient, which we assume to be constant, and $\zeta$ is a Gaussian white noise with
\begin{equation}
\left<\zeta(x,t)\right> = 0 \,
\end{equation}
and
\begin{equation}\label{eq:varwhitenoise1D}
\left<\zeta(x,t) \zeta(x',t')\right> = \delta(x-x')\delta(t-t') \,.
\end{equation}
The results that will be presented here correspond to dimensionless variables, with a domain length $L=1$ and reference time $T=L^2/D=1$.
Periodic boundary conditions
\begin{equation}
u(x=L,t)=u(x=0,t) \, ,
\end{equation}
are prescribed at the ends of the computational domain,  and the initial conditions are defined as
\begin{equation}
u(x,t=0)=u_{0} \, . 
\end{equation}
The periodic boundary conditions guarantee the conservation of mass (or, equivalently, of the number of particles) in the domain, so that the average concentration in the domain is always $u_0$.

If we linearize the equation around a large value of the average concentration $u_0$,  
the equal-time autocorrelation function of the fluctuations at points $x$ and $x'$ for an infinite domain can be expressed as \cite{de2015finite}
\begin{equation}\label{eq:theor_corr_eq}
	\left<u'(x,t) u'(x',t)\right> =   u_{0} \delta(x-x')  \,,
\end{equation}
where $u'=u-\left<u\right>=u-u_0$,
and the so-called static structure factor, which is the Fourier transform of the equal-time correlation function \cite{de2006hydrodynamic}, becomes 
\begin{equation}\label{eq:theor_struc_eq}
	S(k) = \left<\mathcal{U}(k,t) \mathcal{U}(-k,t)\right> =u_{0}  
\end{equation}
where $k$ is the wavenumber, and where the Fourier transform $\mathcal{U}(k)$ of $u(x)$ is defined as
\begin{equation}\label{eq:fourier}
	\mathcal{U}(k) = \int_{-\infty}^{\infty}{u(x)\exp{(-i k x)}} dx  \,.
\end{equation}

%%%%

\subsection{Discrete model, simulation setup and discrete structure factor}
\label{sec:struc}

We apply the numerical approximation described in Section \ref{sec:FE} to Equation \eqref{eq:diffusion}. The one-dimensional domain of length $L$ is split into $N_e$ identical elements. Lagrange polynomials are used as shape functions, both linear ($p1$-elements) and quadratic ($p2$-elements). When $p1$-elements are used, the total number of nodes is $N_n=N_e+1$; when $p2$-elements are used, the number of nodes is $N_n=2 N_e+1$. Since the nodes are always equidistant, the distance between two nodes regardless of the type of element is $\Delta x=L/(N_n-1)$. The number of degrees of freedom is $N_{\text{dof}}\equiv N_n$. 

We use $\widetilde{u}^0_i=u_0$ for every $i=1,..,N_{\text{dof}}$ as the initial condition, where the value of the reference concentration $u_0$ is selected to be high enough for the behavior of the problem to be approximately linear.
We then simulate an initial number of time steps $L^2/(D \Delta t)$, so that the solution reaches an equilibrium state before starting to collect results, and, finally, we run the simulations for a number $N_t$ of time steps, for which data are collected. 

The periodic boundary conditions are implemented by adding Lagrange multipliers to the discrete finite element system \eqref{eq:matricial2}, to force the solution at the extremes of the domain to be identical. These boundary conditions ensure that the total number of particles in the domain remains constant (up to a certain numerical error). This is verified in the computations by integrating the results for the concentration over the domain. For the computations presented in this work, the standard deviation of the total computed number of particles in the domain varies with the discretization parameters, but it is typically smaller or of the order of $ \sim\num{1e-6} \,u_0 L$.

As a consequence of the mass conservation in the domain, the structure factor for wavenumber $k=0$ is zero, unlike the structure factor of the case with an infinite domain given by Equation \eqref{eq:theor_struc_eq}. For the rest of discrete wavenumbers $k_m=2\pi (m-1)/L$, with $m>1$, and provided that the discretization does not introduce artificial correlations in the solution, the structure factor of the solution should approximate that of Equation \eqref{eq:theor_struc_eq}, $S=u_0$. This is equivalent to saying that, although the solution for the continuum equation in an infinite domain has a Gaussian probability distribution, the theoretical {discrete} solution in our bounded one-dimensional domain follows a multinomial distribution corresponding to $N_p=u_0 L$ particles in $N_n-1$ bins. It is worth noting that, as the discretization is refined (i.e. as $N_n$ increases), the multinomial distribution becomes closer to a Gaussian distribution. 

However, the discretization introduces artificial correlations in the solution, so the structure factor of the discrete FE results will not necessarily converge to the theoretical one.
An analysis of the one-dimensional discretized equation with uniform linear ($p1$) Lagrange elements shows instead that the numerical results with $p1$-elements will converge to a solution with a discrete structure factor that can be expressed as \cite{de2015finite}
\begin{equation}\label{eq:theor_struc_eq_FE}
S_{\text{th, FE-p1}}(k_m) =  \frac{9 u_{0}}{\left[2+\cos{(k_m \Delta x)}\right]^2}
\sum_{j\in\mathbb{Z}}{\left(\frac{\sin{\left(\frac{k_m\Delta x}{2} -\pi j\right)}}{\frac{k\Delta x}{2}-\pi j} \right)^4}   \,,
\end{equation}   
where $\mathbb{Z}$ represents the set of integer numbers.
The structure factor given by Equation\eqref{eq:theor_struc_eq_FE} only converges to the continuum structure factor given by Equation \eqref{eq:theor_struc_eq} in the limit $k_m\Delta x \rightarrow 0$.

To evaluate our numerical approach, we compare the theoretical structure factor with the structure factor of the numerical solution, computed as 
\begin{equation}\label{eq:discreteS}
S(k_m)=\left<{\widetilde{\mathcal{U}}}(k_m){{\widetilde{\mathcal{U}}}^*}(k_m)\right>  \,,
\end{equation}
where 
$\widetilde{\mathcal{U}}(k_m)$ is the discrete Fourier transform of the normalized discrete concentration, defined as
\begin{equation}\label{eq:unitarydfft}
{\widetilde{\mathcal{U}}}(k_m) = \frac{1}{\sqrt{L}} \sum_{n=1}^{N_{n}-1}   {\widetilde{u}}_n \Delta \widetilde{L}_n \exp\left(-\frac{i 2 \pi (m-1) (n-1)}{N_{n}-1}\right) \,,
 \end{equation}
 ${\widetilde{\mathcal{U}}}^*(k_m)$ is its complex conjugate and $\Delta \widetilde{L}_n$ is the equivalent volume (equivalent length in 1D) associated to node $n$,
\begin{equation}
\Delta \widetilde{L}_n = \int_{\Omega}{\phi_n(x)}dx\,.
\end{equation}

%%%

\subsection{Results with a semi-implicit time integrator}\label{sec:results_eq_semi}

\begin{figure}
\centering
\includegraphics[width=6.7cm]{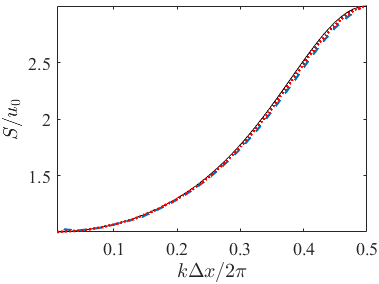}
\includegraphics[width=6.7cm]{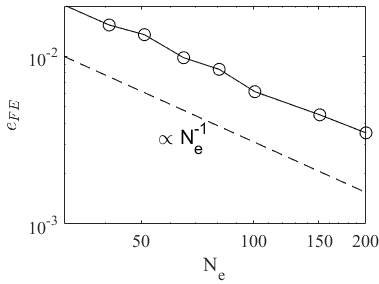}
\caption{Left: Structure factor of the discrete solution; theoretical structure factor for the discrete equation (solid line) and computed with $N_n=51$  (dashed) and $N_n=101$  (dotted). Right: error ${e_{FE}}$ in the computed discrete structure factor vs the number of elements ($\circ$). All computations have $p1$-elements, $\Delta t=\num{1e-4}$, $N_t=\num{1e6}$, $u_{0}=10000$. It can be seen that the solution obtained from the FE computation converges to the theoretical solution of the \emph{discretized} equation as the number of elements increases.}
\label{fig:strfactor_correlated}
\end{figure}

Figure \ref{fig:strfactor_correlated} shows the discrete structure factor based on the finite element results using linear Lagrange elements and a semi-implicit integration scheme.
As explained in Section \ref{sec:timeintegrators}, a semi-implicit integration scheme has the advantage of making the spatial correlations independent from the time step.
As shown in the figure, the discrete structure factor converges towards the theoretical structure factor of the discrete weak equation, given by Equation \eqref{eq:theor_struc_eq_FE}. The error
\begin{equation}
    e_{\text{FE}}=\frac{1}{N_k}\sum_{m=1}^{N_k}{\frac{\left|S(k_m)-S_{\text{FE-p1}}(k_m)\right|}{\left|S_{\text{FE-p1}}(k_m)\right|}} \,,
\end{equation}
where $k_m=2\pi m/L$ (for $m=1,2,..,(N_{n}-1)/2$) are the $N_k=(N_{n}-1)/2$ discrete wavenumbers, 
is reduced with increasing number of elements $N_e$ following a power law ${e_{FE}}\propto N_e^{-1}$. 

\begin{figure}
\centering
\includegraphics[width=13.cm]{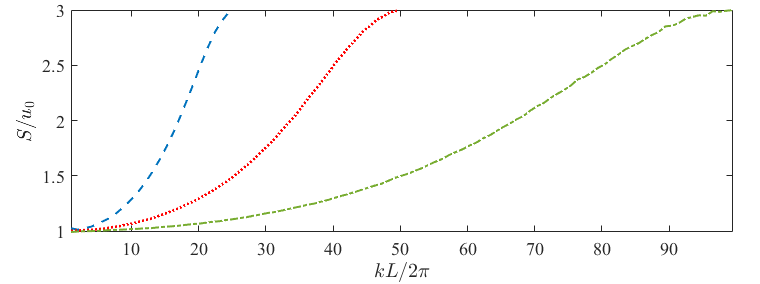}
\caption{Structure factor of the discrete solution as a function of $k L$; computed with $N_n=51$  (dashed), $N_n=101$  (dotted) and $N_n=201$  (dash-dotted). All computations with $p1$-elements, $\Delta t=\num{1e-4}$, $N_t=\num{1e6}$, $u_{0}=10000$.}
\label{fig:strfactor_correlated_vskL}
\end{figure}

However, as explained in Section \ref{sec:struc}, this theoretical discrete structure factor towards which the finite element-solution converges is different from the continuum one, which is a constant ($S/u_0=1$). This is due to the fact that the spatial discretization introduces artificial correlations in the discrete solution; as a result, the structure factor only converges to that of the continuum in the low wavenumber limit. In Figure \ref{fig:strfactor_correlated_vskL}, the same results are represented, but as a function of $k L$ instead of $k \Delta x$. We see that, as the discretization is refined, the structure factor is close to the continuum one for a larger interval of wavenumbers. Nevertheless, artificial correlations are always present in the solution, and the resulting variance cannot be interpreted independently of the discretization.   
It is worth noting that the obtained results agree with the finite element solution obtained by de la Torre and co-workers \cite{de2015finite}, who use a Petrov-Galerkin discretization and also obtain a discrete structure factor that converges to that given by Equation \eqref{eq:theor_struc_eq_FE}. The reason for this is that the artificial correlations only depend on the shape functions used to approximate the solution $u$ in Equation \eqref{eq:discreteu}, which are linear Lagrange polynomials in both works. 

\begin{figure}
\centering
\includegraphics[width=6.7cm]{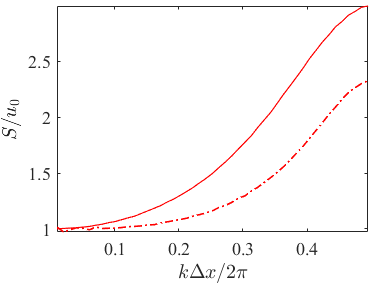}
\caption{Structure factor of the discrete solution computed with $p1$-elements (solid line)  and $p2$-elements (dash-dot). The computations have been run with $N_{n}=101$, $\Delta t=\num{1e-4}$, $N_t=\num{1e6}$, $u_{0}=10000$.}
\label{fig:strfactor_correlated_p1p2}
\end{figure}

\begin{figure}
\centering
\includegraphics[width=6.7cm]{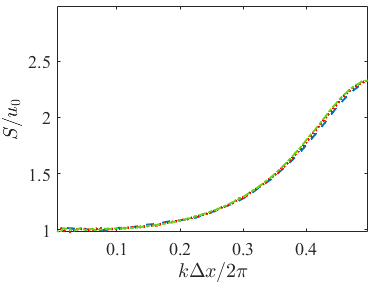}
\includegraphics[width=6.7cm]{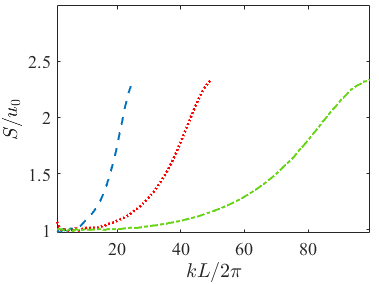}
\caption{Structure factor of the discrete solution obtained with $p2$-elements, computed with $N_{n}=51$  (dashed) and $N_{n}=101$  (dotted). The computations have been run with $\Delta t=\num{1e-4}$, $N_t=\num{1e6}$, $u_{0}=10000$.}
\label{fig:strfactor_correlated_p2}
\end{figure}

Figure \ref{fig:strfactor_correlated_p1p2} shows results obtained with both linear and quadratic elements. It can be seen that increasing the order contributes to approximating better the structure factor towards the continuum one. Although the convergence is still limited to the low wavenumber limit, the results obtained with quadratic elements present a larger interval of wavenumbers for which the error in the structure factor is below a given threshold.
Figure \ref{fig:strfactor_correlated_p2} shows the results obtained with quadratic elements and different mesh sizes. Like the results obtained with linear elements, the results computed with quadratic elements do not converge towards the continuum as the mesh is refined, indicating the presence of artificial correlations induced by the discretization in the solution.

 {
The results presented in this section have been focused on the static structure factor because it serves the purpose of showing the presence of artificial spatial correlations (or eventual lack thereof) in the solution, which is the focus of this work. Nevertheless, in \ref{sec:dynamicSF}, complementary results are presented showing the dynamic structure factor, which contains information on the time evolution of the solution as well.
}

\subsection{Results with fully implicit and fully explicit time integrators}\label{sec:results_eq_impli}

\begin{figure}
\centering
\includegraphics[width=6.7cm]{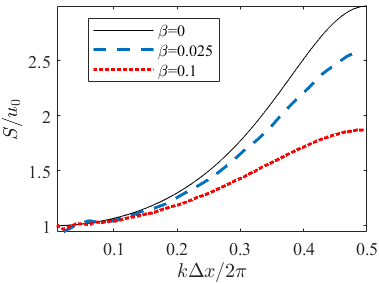}
\includegraphics[width=6.7cm]{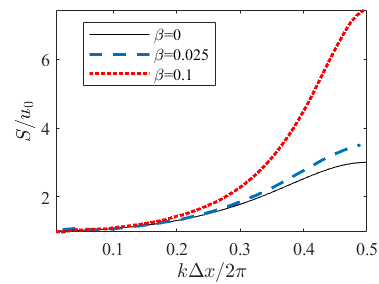}
\caption{Structure factor of the discrete solution with $p1$-elements, computed with $N_n=51$  (dashed) and $N_n=101$  (dotted). Left: fully implicit time integration ($\alpha=0$). Right: fully explicit time integration ($\alpha=1$). All computations with $\Delta t=\num{1e-5}$, $N_t=\num{1e6}$, $u_{0}=10000$. 
As a reference, numerical results obtained with $\alpha=1/2$ are also displayed (solid line), which are equivalent to the case with $\beta=0$.}
\label{fig:strfactor_correlated_beta}
\end{figure}

\begin{figure}
\centering
\includegraphics[width=6.7cm]{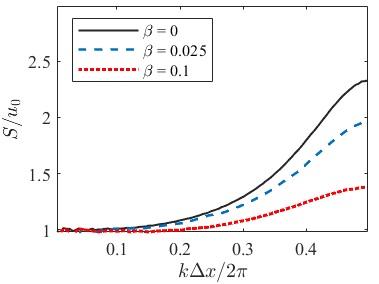}
\includegraphics[width=6.7cm]{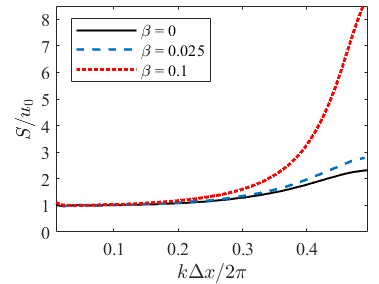}
\caption{Structure factor of the discrete solution with $p2$-elements, computed with $N_{n}=51$  (dashed) and $N_{n}=101$  (dotted). Left: fully implicit time integration ($\alpha=0$). Right: fully explicit time integration ($\alpha=1$). All computations with $\Delta t=\num{1e-5}$, $N_t=\num{1e6}$, $u_{0}=10000$. As a reference, numerical results obtained with $\alpha=1/2$ are also displayed (solid line), which are equivalent to the case with $\beta=0$.}
\label{fig:strfactor_correlated_beta_p2}
\end{figure}

The results shown so far for the semi-implicit time integrator present artificial correlations that are related to the spatial discretization. However, as explained in Section \ref{sec:timeintegrators}, if instead of the semi-implicit integration scheme, one uses a fully explicit or a fully implicit scheme, additional spatial correlations dependent on the dimensionless time step $\beta=D \Delta t/\Delta x^2$ appear in the solution \cite{donev2010accuracy}. 
Figures \ref{fig:strfactor_correlated_beta}
and \ref{fig:strfactor_correlated_beta_p2} show the results obtained for linear and quadratic elements, using fully explicit and fully implicit integration schemes. The results confirm that indeed the obtained structure factor is influenced by additional correlations that depend on the dimensionless time step $\beta$. It is worth noting that the influence of the time step on the stability is an independent issue altogether: the implicit method is unconditionally stable, while the explicit becomes unstable if the time step exceeds a certain threshold.

The results obtained with the fully explicit and fully implicit integration schemes are consistent with what other authors have shown before \cite{donev2010accuracy, delong2013temporal}, and illustrate the advantage of the semi-implicit integration scheme, as it does not require decreasing the time step to yield the correct spatial structure. In the rest of this work, only results obtained with the semi-implicit scheme ($\alpha=1/2$) are presented.

%%%

\section{A linear mapping to remove artificial correlations} \label{sec:decorrelation}

As explained in Sections \ref{sec:timeintegrators} and \ref{sec:coarsegraining}, the spatial discretization can introduce artificial correlations in the numerical solution. These artificial correlations do not represent a numerical error; they are instead the natural consequence of the choice of shape functions for the spatial discretization, and they are thus part of a consistent solution in the context of that basis. In spite of this, these discretization-related correlations make the discrete solution difficult to interprete, and can also make non-linear terms in fluctuating-hydrodynamics equations difficult to handle.
In this section, we propose a linear transformation between the numerical solution and an equivalent discrete solution free of these artificial correlations. The goal is to perform a change of basis to obtain an equivalent discrete solution $\widetilde{\widetilde{\mathbf{u}}}$ that approximates the solution as
\begin{equation}\label{eq:discreteu2}
u(\mathbf{x}, t) \approx \sum_{i=1}^{N_{\text{dof}}}{\widetilde{\widetilde{u}}_i(t) \psi_i(\mathbf{x})} \,,
\end{equation} 
where $\psi_i$ are basis functions defined so as not to introduce artificial spatial correlations in the finite element solution, that is, so as to satisfy
\begin{equation}\label{eq:uncorrelated_weak}
    \int_{\Omega} \psi_i(\mathbf{x}) \psi_j(\mathbf{x}) d^3\mathbf{x}= 0 \,, \,\,\,\, \text{if } i\neq j \,.
\end{equation}

\subsection{Spatial decorrelation matrix}
\label{sec:decorrelationMatrix}

In this section, we show that a linear relationship exists between the discrete solution $\widetilde{\mathbf{u}}$ and an equivalent discrete solution $\widetilde{\widetilde{\mathbf{u}}}$ that is free of the artificial correlations introduced by the spatial discretization. 
We start by re-writing here Equation \eqref{eq:discreteu}
\begin{equation*}
u(\mathbf{x}, t) \approx \sum_{i=1}^{N_{\text{dof}}}{\widetilde{u}_i(t) \phi_i(\mathbf{x})} \,,
\end{equation*}
where $\phi_i(\mathbf{x})$ are the basis functions used to approximate the solution, and 
Equation \eqref{eq:discreteu2}, which provides an alternative approximation using shape functions $\psi_i$, 
\begin{equation*}
u(\mathbf{x}, t) \approx \sum_{i=1}^{N_{\text{dof}}}{\widetilde{\widetilde{u}}_i(t) \psi_i(\mathbf{x})} \,.
\end{equation*}
In order to satisfy Equation \eqref{eq:uncorrelated_weak} the shape functions $\psi_i(\mathbf{x})$ are defined such that matrix $\mathbf{M}_{\psi \psi}$ 
\begin{equation}\label{eq:Ldecorr0}
({M}_{\psi \psi})_{ij}=\int_{\Omega}{\psi_i \psi_j}d^3\mathbf{x} \, 
\end{equation}
is diagonal, with
\begin{equation}\label{eq:Ldecorr}
({M}_{\psi \psi})_{ii}=\int_{\Omega}{\psi_i \psi_i}d^3\mathbf{x}=\Delta \widetilde{\widetilde{V}}_i \, ,
\end{equation}
where $\Delta \widetilde{\widetilde{V}}_i$ is an equivalent discrete volume associated with each degree of freedom $i$,
\begin{equation}\label{eq:voldeci}
\Delta \widetilde{\widetilde{V}}_i = \int_{\Omega}{\psi_i}d^3\mathbf{x}\,.
\end{equation}
Note that we have not defined the shape functions $\psi_i$ beyond saying that they lead to a diagonal mass matrix, so $\Delta \widetilde{ \widetilde{ {V}_i}}$ is undefined so far. We will discuss later in this section which definition of $\Delta \widetilde{\widetilde{{V}_i}}$ is most advantageous (see Equation \eqref{eq:yisq}).

Since the discrete solutions ${\widetilde{\mathbf{u}}}$ and $\widetilde{\widetilde{\mathbf{u}}}$ both approximate the same continuum function, we require them to satisfy 
\begin{equation}
\sum_{j=1}^{N_{\text{dof}}}{{\widetilde{u}}_j(t) {{\phi}}_j(\mathbf{x})}
=
\sum_{j=1}^{N_{\text{dof}}}{\widetilde{\widetilde{u}}_j(t) \psi_j(\mathbf{x})} \,.
\end{equation}
Multiplying by $\phi_i$ and integrating over the domain leads to
\begin{equation}
\int_{\Omega}{{{\phi}}_i \sum_{j=1}^{N_{\text{dof}}} {{\phi}}_j} \widetilde{u}_j d^3\mathbf{x}
=
\int_{\Omega}{{{\phi}}_i \sum_{j=1}^{N_{\text{dof}}} \psi_j} \widetilde{\widetilde{u}}_j d^3\mathbf{x} \,,
\end{equation}
which is equivalent to the following equation in matrix form 
\begin{equation}\label{eq:key}
\mathbf{A}_{\phi}^T \widetilde{\mathbf{u}} = \mathbf{A}_{\psi} \, \widetilde{\widetilde{\mathbf{u}}} \,,
\end{equation}
where $\mathbf{A}_{\psi}$ is a diagonal matrix satisfying 
\begin{equation}
\mathbf{M}_{\psi \psi}=\mathbf{A}_{\psi}\,\mathbf{A}_{\psi} \,,
\end{equation}
and $\mathbf{A}_{\phi}^T$ is the transpose of $\mathbf{A}_{\phi}$, satisfying
\begin{equation}\label{eq:Ldecomp}
\mathbf{M}_{\phi \phi}=\mathbf{A}_{\phi}\,\mathbf{A}_{\phi}^T \,,
\end{equation}
with   matrix $\mathbf{M}_{\phi \phi}$ defined as
\begin{equation}\label{eq:Lcorr0}
({M}_{\phi \phi})_{ij}=\int_{\Omega}{\phi_i \phi_j}d^3\mathbf{x} \,. 
\end{equation}
It is worth pointing out that, since in Section \ref{sec:FE} we have used a standard Galerkin discretization with the same shape functions for test function and solution, $\mathbf{M}_{\phi \phi}$ and the mass matrix $\mathbf{M}$ of the finite element system are the same in our numerical setup. However, this is not necessarily the case for a general discretization, in which different shape functions may be used to approximate the test function. Therefore, $\mathbf{M}_{\phi \phi}$ is not necessarily equal to the mass matrix of the finite element system leading to the solution $\widetilde{\mathbf{u}}$.

Matrix $\mathbf{A}_{\phi}$, resulting from the decomposition given by Equation \eqref{eq:Ldecomp}, can be expressed as 
\begin{equation}\label{eq:AL}
\mathbf{A}_{\phi} = \mathbf{U}_{\phi} \sqrt{\mathbf{\Sigma}_{\phi}} \mathbf{U}_{\phi}^T \mathbf{U}_*^T\,,
\end{equation} 
where $\mathbf{U}_*$ and $\mathbf{U}_{\phi}$ are unitary matrices and $\mathbf{\Sigma}_{\phi}$ is a diagonal matrix satisfying
\begin{equation}
\mathbf{M}_{\phi \phi} = \mathbf{U}_{\phi} \mathbf{\Sigma}_{\phi} \mathbf{U}_{\phi}^T\,.
\end{equation}
Equation \eqref{eq:Ldecomp} holds for any unitary matrix $\mathbf{U}_*$; the decomposition \eqref{eq:Ldecomp} is therefore not unique. However, not every unitary matrix $\mathbf{U}_*$ will in addition guarantee that the mass is conserved.
Mass conservation requires
\begin{equation}
\sum_{j=1}^{N_{\text{dof}}} \int_{\Omega}{ {{\phi}}_j} \widetilde{u}_j d^3\mathbf{x}
=
\sum_{j=1}^{N_{\text{dof}}} \int_{\Omega}{\psi_j} \widetilde{\widetilde{u}}_j d^3\mathbf{x} \,,
\end{equation}
which in matrix form can be written as
\begin{equation}\label{eq:massconservationmatricial}
\Delta{\widetilde{\mathbf{ V}}}^T\,\widetilde{\mathbf{u}} 
=
\Delta \widetilde{\widetilde{\mathbf{ V}}}^T\,\widetilde{\widetilde{\mathbf{u}}} \,,
\end{equation}
where $\Delta \widetilde{\widetilde{\mathbf{ V}}}$ and $\Delta {\widetilde{\mathbf{ V}}}$ are vectors with the corresponding equivalent volume associated with each degree of freedom, with $\Delta \widetilde{\widetilde{\mathbf{ V}}}$ defined by Equation \eqref{eq:voldeci} and $\Delta {\widetilde{\mathbf{ V}}}$ defined by
\begin{equation}\label{eq:voli}
\Delta \widetilde{ V}_i
=
\int_{\Omega}{{{\phi}}_i} d^3\mathbf{x} \,.
\end{equation}
If we now use Equations \eqref{eq:key} and \eqref{eq:AL} inside Equation \eqref{eq:massconservationmatricial}, we obtain
\begin{equation}
\Delta{\widetilde{\mathbf{ V}}}^T\,\widetilde{\mathbf{u}} 
=
\Delta \widetilde{\widetilde{\mathbf{ V}}}^T\,\mathbf{A}_{\psi}^{-1}
\, \mathbf{U}_*\,\mathbf{U}_{\phi}\,\sqrt{\mathbf{\Sigma}_{\phi}}\,\mathbf{U}_{\phi}^T \,\widetilde{\mathbf{u}} \,.
\end{equation}
Given that this identity must be satisfied for any vector  $\widetilde{\mathbf{u}}$, we obtain
\begin{equation}
\mathbf{U}_* \,\widetilde{\mathbf{y}} 
= \widetilde{\widetilde{\mathbf{y}}} \,, 
\end{equation}
where ${\widetilde{\mathbf{y}}}$ and $\widetilde{\widetilde{\mathbf{y}}}$ are vectors defined as
\begin{equation}
{\widetilde{\mathbf{y}}}=  \mathbf{U}_{\phi}\,\sqrt{\mathbf{\Sigma}^{-1}_{\phi}} \,\mathbf{U}_{\phi}^T \, \Delta \widetilde{\mathbf{V}}\,, 
\end{equation}
and 
\begin{equation}
\widetilde{\widetilde{\mathbf{y}}}=  \mathbf{A}_{\psi}^{-1}\, \Delta \widetilde{\widetilde{\mathbf{V}}}\,. 
\end{equation}
Therefore, matrix $\mathbf{U}_*$ can be defined as the 
rotation matrix that transforms vector ${\widetilde{\mathbf{y}}}$ into $\widetilde{\widetilde{\mathbf{y}}}$.

With these definitions, we can rewrite Equation \eqref{eq:key} as
\begin{equation}\label{eq:key2}
\widetilde{\widetilde{\mathbf{u}}} = \mathbf{A}_{\psi}^{-1} \,\mathbf{A}_{\phi}^T \,\widetilde{\mathbf{u}} = \mathbf{A}_{\psi}^{-1} \,\mathbf{U}_*\,\mathbf{U}_{\phi}\,\sqrt{\mathbf{\Sigma}_{\phi}}\,\mathbf{U}_{\phi}^T \,\widetilde{\mathbf{u}}   \,,
\end{equation}
which provides a mass-preserving linear relationship between the discrete solution $\widetilde{\mathbf{u}}$ and a discrete solution $\widetilde{\widetilde{\mathbf{u}}}$ that  does not have artificial spatial correlations introduced by the spatial discretization:
\begin{equation}\label{eq:mappingQ}
   \widetilde{\widetilde{\mathbf{u}}} = \mathbf{Q}\,\widetilde{\mathbf{u}}\,,
\end{equation}
where matrix $\mathbf{Q}$ is given by
\begin{equation}
\mathbf{Q} = \mathbf{A}_{\psi}^{-1} \,\mathbf{A}_{\phi}^T
= \mathbf{A}_{\psi}^{-1} \,\mathbf{U}_*\,\mathbf{U}_{\phi}\,\sqrt{\mathbf{\Sigma}_{\phi}}\,\mathbf{U}_{\phi}^T    \,
\end{equation}
and can be interpreted as a decorrelation matrix that removes the artificial spatial correlations from the discrete solution. This decorrelation matrix can be applied as a postprocessing step to remove the artificial spatial correlations from the numerical solution. 

As mentioned above,  $\Delta \widetilde{ \widetilde{ \mathbf{V}}}$ is undefined so far.  A convenient definition is given by
\begin{equation}\label{eq:yisq}
    \Delta \widetilde{ \widetilde{ \mathbf{V}}}_i
    \equiv \widetilde{\mathbf{y}}_i^2 \,.
\end{equation}
This definition ensures that $\widetilde{\mathbf{y}}$ and $\widetilde{\widetilde{\mathbf{y}}}$ are the same, and, therefore, $\mathbf{U}_*$ becomes the identity matrix. As a result, $\mathbf{A}_{\phi}$ becomes symmetric, and the decorrelation matrix becomes 
\begin{equation}\label{eq:mappingD}
\mathbf{Q} = \mathbf{A}_{\psi}^{-1} \,\mathbf{A}_{\phi}^T
= \mathbf{A}_{\psi}^{-1} \,\mathbf{U}_{\phi}\,\sqrt{\mathbf{\Sigma}_{\phi}}\,\mathbf{U}_{\phi}^T    \,.
\end{equation}
An advantage of this definition is that the decorrelation matrix can be approximated by a sparse matrix with a higher sparsity degree than what is obtained when the rotation matrix $\mathbf{U}_*$ is a dense matrix. In Section \ref{sec:computaspects}, we discuss more details about the implementation and computational aspects of the mapping given by Equation \eqref{eq:mappingD}. 

So far, we have not considered any effect of the boundary conditions on the decorrelation matrix, and, in general, the analysis presented in this section is valid for any boundary condition. However, in some cases, it may be desirable to use a decorrelation matrix based on a matrix $\mathbf{M}_{\phi \phi}$ and vector $\Delta \widetilde{\mathbf{V}}$ that have been modified to incorporate the effect of the boundary conditions. For instance, since here we are solving problems with periodic boundary conditions, in order to maintain the periodicity of the discrete solution, only $N_{\text{dof}}-1$ elements of $\widetilde{\mathbf{u}}$ are mapped, thereby excluding the solution of one of the extremes of the domain.
Matrix $\mathbf{M}_{\phi \phi}$ and vector $\Delta \widetilde{\mathbf{V}}$ become, respectively, an $(N_{\text{dof}}-1)\times (N_{\text{dof}}-1)$ matrix and a vector of length $(N_{\text{dof}}-1)$, and the elements in each eliminated row or column are added in suitable positions of the row or column corresponding to the other extreme of the domain, so as to reflect the periodicity of the domain in the arrays. 

It is worth noting that the linear mapping given by Equation \eqref{eq:mappingQ} depends only on the shape functions used to approximate the solution. This means that the choice of finite-element discretization does not need to be based on whether it produces the correct structure factor, because, once a discrete solution is obtained, an equivalent solution with a correct structure factor can be found through the proposed linear transformation.  
Moreover, this makes explicit that the coarse-graining and the spatial discretization used to solve the diffusion equation numerically are two separate things, and that the spatial discretization should therefore be designed on the basis of numerical analysis exclusively. 

In Section \ref{sec:results_diff2_decorr}, we present results for the one-dimensional boundary value problem discussed in Section \ref{sec:results_diff2}, which are obtained by applying the linear mapping presented in this section to remove the artificial correlations introduced by the spatial discretization.

%%%

\subsection{Implementation and computational aspects}
\label{sec:computaspects}

 {The decorrelation matrix proposed in Section \ref{sec:decorrelationMatrix} is in general a dense matrix. Although this is not an issue for the 1D problems solved in this work, the mapping can lead to a high computational cost for models with a large number of degrees of freedom.} To preserve the computational advantages of the finite element method, the matrix should be sparse.

\begin{figure}
\centering
\includegraphics[width=6.7cm]{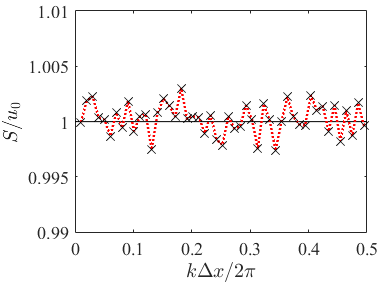}
\includegraphics[width=6.7cm]{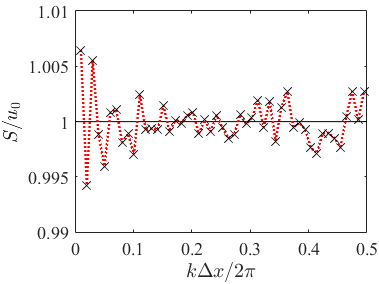}
\caption{Structure factor of the discrete solution after mapping with the decorrelation matrix; computed with a dense decorrelation matrix (dotted) and with a sparse decorrelated matrix with threshold $\epsilon=\num{1e-5}$  ($\times$). Left: $p1$-elements; right: $p2$-elements. All computations with semi-implicit time integration, $N_{\text{dof}}=101$, $\Delta t=\num{1e-4}$, $N_t=\num{1e6}$, $u_{0}=10000$.}
\label{fig:strfactor_truncated}
\end{figure}

If we use the definition of the decorrelation matrix $\mathbf{Q}$ given by Equation \eqref{eq:mappingD}, with the equivalent volumes $\Delta \widetilde{\widetilde{V}}_i$ defined by Equation \eqref{eq:yisq}, the matrix can be approximated by a sparse matrix, by setting to zero the elements of the matrix that have an absolute value below a threshold $\epsilon$.  {This is a reasonable approximation in our case, because the matrix $\mathbf{Q}$ that we obtain already has many elements with small value, corresponding to degrees of freedom belonging to two different elements that are not contiguous and that are relatively far from each other. It is worth pointing out that the sparse matrix that we obtain in this way is still a full-rank matrix.}

Figure \ref{fig:strfactor_truncated} shows the structure factor obtained with such truncated matrices with a threshold $\epsilon=\num{1e-5}$. The truncated matrix gives results that are very close to those of the original dense matrix (with a maximum added relative error in the computed structure factor around $\num{5e-5}$, which is significantly smaller than the typical numerical error in the solution).

\begin{figure}
\centering
\includegraphics[width=6.7cm]{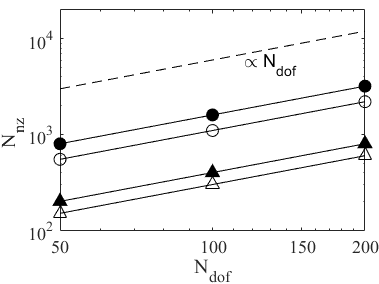}
\caption{Number $N_{\text{nz}}$ of non-zero elements as a function of the number of degrees of freedom $N_{\text{dof}}$ for mass matrix $\mathbf{M}$ with $p1$-elements ($\triangle$) and $p2$-elements ($\blacktriangle$), truncated decorrelation matrix $\mathbf{Q}$ ($\epsilon=\num{1e-5}$) with $p1$-elements ($\circ$) and $p2$-elements ($\bullet$).}
\label{fig:computComplexity}
\end{figure}

Figure \ref{fig:computComplexity} shows the number of non-zero elements in the decorrelation matrix with both $p1$- and $p2$-elements, computed as well with threshold $\epsilon=\num{1e-5}$. Although the number of non-zero elements in the decorrelation matrices is higher than in the respective mass matrices, they both increase linearly with the number of degrees-of-freedom $N_{\text{dof}}$. This indicates that the computational complexity of the problem is not increased by the mapping, thus preserving the efficiency of the computational method.

An example of implementation of a finite-element method including a decorrelation step for the problems discussed in previous sections consists of the following steps:
\begin{enumerate}
    \item Discretize the domain: define the elements and node connectivity, select the shape functions
    \item Assemble system matrices $\mathbf{M}$ and $\mathbf{D}$ 
    \item Implement the components of the boundary conditions (e.g. add Lagrange multipliers for periodic boundary conditions) and source terms (e.g. compute matrix $\mathbf{A_D}$ for the fluctuating source term) that are not time-dependent
    \item Compute the decorrelation matrix $\mathbf{Q}$:
    \begin{enumerate}
        \item Compute $\mathbf{U}_{\phi}$ and $\mathbf{\Sigma}_{\phi}$ through a singular value decomposition of $\mathbf{M}_{\phi\phi}$ (note that, since $\mathbf{M}_{\phi\phi}$ is symmetric and positive definite, the eigenvalues and singular values are identical)
        \item Compute $\Delta \widetilde{ \widetilde{ \mathbf{V}}}$ using Equation \eqref{eq:yisq}
        \item Compute $\mathbf{A}_{\psi}$ as a diagonal matrix with diagonal elements given by $\sqrt{\Delta \widetilde{ \widetilde{ \mathbf{V}}}_i}$
        \item Compute the decorrelation matrix as $\mathbf{Q}=\mathbf{A}_{\psi}^{-1}\,\mathbf{U}_{\phi}\,\sqrt{\mathbf{\Sigma}_{\phi}}\,\mathbf{U}_{\phi}^T$
        \item Make $\mathbf{Q}$ sparse by setting elements with absolute value below a given tolerance to zero.
    \end{enumerate}
    \item Initialize the solution $\widetilde{\mathbf{u}}^0$
    \item Start the time-dependent loop, for $n=1$ to $n=N_t$ do:
     \begin{enumerate}
        \item Compute the fluctuating source term
        \item Solve the linear system to obtain the finite-element solution array $\widetilde{\mathbf{u}}^{n}$ as a function of $\widetilde{\mathbf{u}}^{n-1}$ (see Section \ref{sec:timeintegrators})
        \item Compute the mapped solution as $\widetilde{\widetilde{\mathbf{u}}}^{n} = \mathbf{Q}\,\widetilde{\mathbf{u}}^{n}$ and store it
    \end{enumerate}
\end{enumerate}

%%%

\subsection{Results using the spatial decorrelation matrix}
\label{sec:results_diff2_decorr}

\begin{figure}
\centering
\includegraphics[width=6.7cm]{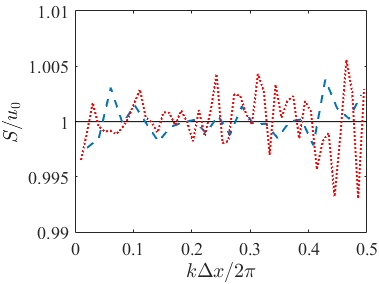}
\includegraphics[width=6.7cm]{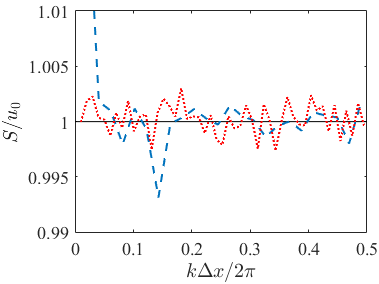}
\includegraphics[width=6.7cm]{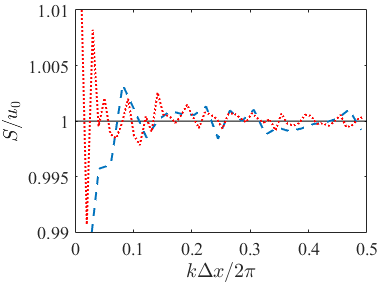}
\includegraphics[width=6.7cm]{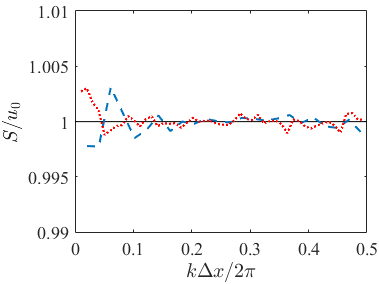}
\caption{Structure factor of the discrete solution after mapping with the decorrelation matrix computed with $N_n=51$  (dashed) and $N_n=101$  (dotted), and theoretical structure factor for the continuum equation (solid line). Top left: $\Delta t=\num{1e-3}$, $N_t=\num{1e6}$; top right: $\Delta t=\num{1e-4}$, $N_t=\num{1e6}$; bottom left: $\Delta t=\num{1e-5}$, $N_t=\num{1e7}$; bottom right: $\Delta t=\num{1e-4}$, $N_t=\num{1e7}$. All computations with $p1$-elements, semi-implicit time integration,  $u_{0}=10000$.}
\label{fig:strfactor_decorr1}
\end{figure}
\begin{figure}
\centering
\includegraphics[width=6.7cm]{figures1/S_decorrMatrix_semi_Ne100_Nt1e6_dt1en4.png}
\includegraphics[width=6.7cm]{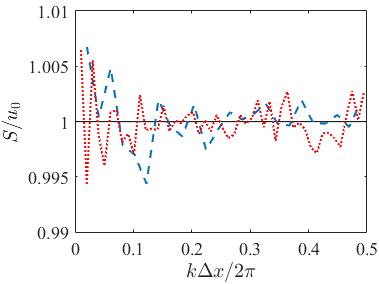}
\caption{Structure factor of the discrete solution after mapping with the decorrelation matrix computed with $N_n=51$  (dashed) and $N_n=101$  (dotted), and theoretical structure factor for the continuum equation (solid line). Left: $p1$-elements; right: $p2$-elements. All computations with semi-implicit time integration, $\Delta t=\num{1e-4}$, $N_t=\num{1e6}$, $u_{0}=10000$.}
\label{fig:strfactor_decorr2}
\end{figure}

Figures \eqref{fig:strfactor_decorr1}
and \eqref{fig:strfactor_decorr2} show the results obtained applying the linear mapping \eqref{eq:key2} given by Equation \eqref{eq:mappingQ} to the finite element results presented in Section \ref{sec:results_diff2}. Only the results obtained with the semi-implicit integration scheme have been used. The discrete structure factor for the decorrelated results is computed as
\begin{equation}
S(k_m)=\left<{\widetilde{\widetilde{\mathcal{U}}}}(k_m){\widetilde{\widetilde{\mathcal{U}}}^*}(k_m)\right>  \,,
\end{equation}
with 
\begin{equation}\label{eq:unitarydfft2}
{\widetilde{\widetilde{\mathcal{U}}}}(k_m) = \frac{1}{\sqrt{L}} \sum_{n=1}^{N_{n}-1}   \widetilde{\widetilde{u}}_n \Delta \widetilde{\widetilde{L}}_n \exp\left(-\frac{i 2 \pi (m-1) (n-1)}{N_{n}-1}\right) \,,
 \end{equation}
where $\Delta \widetilde{\widetilde{L}}_n$ is the equivalent volume (equivalent length in 1D)  related to shape function $\psi$ for node $n$, 
which can be computed using Equation \eqref{eq:yisq}.

The results presented in
Figures \eqref{fig:strfactor_decorr1} have been computed for two different meshes (number of nodes $N_n=51$ and $N_n=101$) with linear shape functions, and for different values of the time step $\Delta t$ and of the number of time steps $N_t$ for which solution data are collected. It can be seen that the linear mapping succeeds in yielding a structure factor that closely approximates that of the continuum ($S/u_0=1$). Once the artificial correlations are removed, and since the solution is completely decorrelated in space, refining the mesh does not improve the solution significantly, and, for the largest time step $\Delta t=\num{1e-3}$ (see top left plot), refining the mesh actually increases the error, since the dimensionless time step $\beta=\Delta t D/\Delta x^2$ becomes larger. For small-enough values of the time step, increasing the number of time steps $N_t$ and the total time of the simulations $t_{\text{tot}}=N_t \Delta t$ is what contributes most to decreasing the error.  

Figure \eqref{fig:strfactor_decorr2} shows results obtained with both linear shape functions ($p1$-elements) and quadratic shape functions ($p2$-elements). Just like refining the mesh does not decrease the numerical error in this case, since the solution is white noise, increasing the order does not influence the accuracy of the results significantly.

 {
In \ref{sec:dynamicSF}, complementary results for the computed dynamic structure structure of the mapped solution are shown, which illustrate that the method not only succeeds in removing artificial spatial correlations from the solution, but is also able to capture the time evolution of the solution.
}
%%%%
%%%

\section{FE solution of a fourth-order 1D diffusion problem}\label{sec:results_diff4}

Following \cite{de2015finite},
we consider next a stochastic diffusion equation of the form given by Equation \eqref{eq:PDEintro} with 
\begin{equation} 
F(u)=u +\left(\frac{\ell_0}{2\pi}\right)^2 \nabla^2 u\,,
\end{equation} which leads to the fourth-order equation
\begin{equation}\label{eq:diff4}
	\frac{\partial {u}}{\partial t} - D \nabla^2 \left(u +\left(\frac{\ell_0}{2\pi}\right)^2 \nabla^2 u\right) = \nabla \cdot \left( \sqrt{2 D u}\boldsymbol{\zeta} \right)\,.
\end{equation}
Unlike Equation \eqref{eq:diff2}, the solution of which has only spatially decorrelated fluctuations, Equation \eqref{eq:diff4} has a solution that exhibits spatial correlations with a finite correlation length defined by $\ell_0$. This allows us to investigate the effect of the mapping for different ratios of cell size to correlation length ${\Delta x}/{\ell_0}$.

\subsection{One-dimensional boundary value problem}
We consider the one-dimensional version of Equation \eqref{eq:diff4}
\begin{equation}\label{eq:diffusion4}
	\frac{\partial u({x},t)}{\partial t} -  D \frac{\partial^2}{\partial x^2}\left(u({x},t) + \left(\frac{\ell_0}{2\pi}\right)^2 \frac{\partial^2 u({x},t)}{\partial x^2}\right)   = \frac{\partial}{\partial x} \left( \sqrt{2 D u_0} {\zeta}({x},t)  \right)\,,
\end{equation}
where $u_0$ and $\ell_0$ are constants. The rest of the variables are defined in the same way as for the second-order problem discussed in Section \ref{sec:results_diff2}.   

Instead of solving Eq. \eqref{eq:diffusion4} directly, we replace the fourth-order PDE with the coupled set of second-order equations
\begin{subequations}
\begin{align}
\frac{\partial u}{\partial t} -  D \frac{\partial^2}{\partial x^2}\left(u + \left(\frac{\ell_0}{2\pi}\right)^2 w\right)   &= \frac{\partial}{\partial x} \left( \sqrt{2 D u_0} \mathbf{\boldsymbol{\zeta}} \right)
\\
w - \frac{\partial^2 u}{\partial x^2} &= 0 \,.
\end{align}
\end{subequations}
This allows us to use the same kind of discretization described in Section \ref{sec:FE}.
Periodic boundary conditions are prescribed at the boundaries of the domain,
\begin{equation}
u(x=L,t)=u(x=0,t) \, ,    
\end{equation}
\begin{equation}
w(x=L,t)=w(x=0,t) \, ,
\end{equation}
and the initial conditions are defined as
\begin{equation}
u(x,t=0)=u_{0} \, ,
\end{equation}
\begin{equation}
w(x,t=0)=0 \, .
\end{equation}
The theoretical static structure factor of the solution to this problem is given by 
\begin{equation}\label{eq:theor_struc_eq2}
	S(k) = u_{0} \frac{1}{1+\frac{k^2}{k_0^2}}  \,,
\end{equation}
with the cutoff wavenumber defined as 
\begin{equation}
k_0=\frac{2\pi}{\ell_0}\,.
\end{equation}
Unlike the case presented in Sec. \ref{sec:results_diff2}, the covariance of the fluctuating solution to this equation for points $x_1$ and $x_2$ is not proportional to the delta function, but to an exponentially decaying function $\propto u_0 k_0 e^{-k_0 (x_1-x_2)}$ with characteristic length $\ell_0$ \cite{de2015finite}.

This problem has a limit case when $\Delta x/\ell_0 \gg 1$. In this case, the correlation length of the fluctuations is much smaller than the discretization length; this corresponds to the spatially decorrelated case described in Section \ref{sec:results_diff2}, which needs the mapping proposed in Section \ref{sec:decorrelation} to remove the artificial correlations introduced by the discretization. 
The other limit case occurs when $\Delta x/\ell_0 \rightarrow 0$. In this case, the discretization captures the correlation length of the fluctuations completely, the variance of the fluctuations not resolved by the mesh tends to zero, and, therefore, the direct FE solution without mapping converges to the continuum solution.
In the following section, we discuss the results for different ratios  $\Delta x/\ell_0$.

Like in Section \ref{sec:results_diff2}, all the variables are made dimensionless variables, with the reference length $L=1$ and the reference time $T=L^2/D=1$. 

%%%%
\subsection{Discretization and results}

\begin{figure}
\centering
\includegraphics[width=6.7cm]{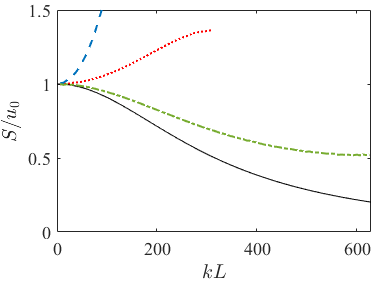}
\includegraphics[width=6.7cm]{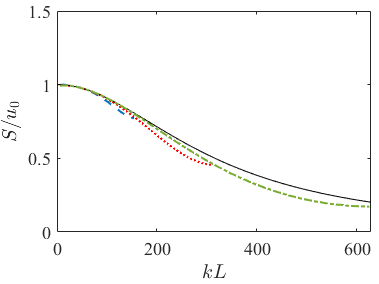}
\caption{Structure factor of the FE results before (left) and after applying the decorrelation matrix (right): theoretical structure factor (solid line), $\Delta x/\ell_0=0.25$ (dash-dot), $\Delta x/\ell_0=0.5$ (dotted) and $\Delta x/\ell_0=1$ (dashed). The results have been computed with $p1$-elements, $\Delta t=\num{1e-4}$, $N_t=\num{5e6}$, $u_0=10000$. }
\label{fig:strfactor_ddft}
\end{figure}

\begin{figure}
\centering
\includegraphics[width=6.7cm]{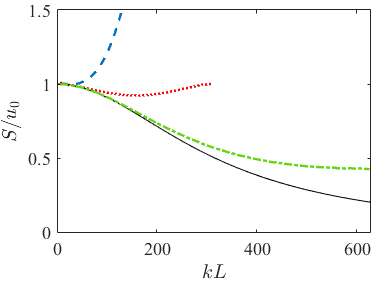}
\includegraphics[width=6.7cm]{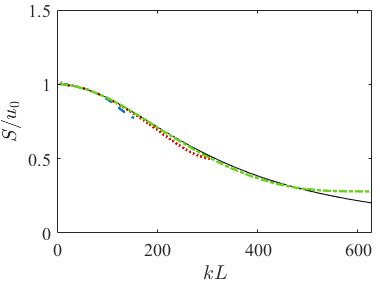}
\caption{Structure factor of the FE results before (left) and after applying the decorrelation matrix (right):  theoretical structure factor (solid line), $\Delta x/\ell_0=0.25$ (dash-dot), $\Delta x/\ell_0=0.5$ (dotted) and $\Delta x/\ell_0=1$ (dashed). The results have been computed with $p2$-elements, $\Delta t=\num{1e-4}$, $N_t=\num{5e6}$, $u_0=10000$. }
\label{fig:strfactor_ddft2}
\end{figure}

The simulation setup for this problem is based on the same discretization schemes and numerical procedures described in Sections \ref{sec:FE} and \ref{sec:struc} for the second-order boundary value problem.
Figures \ref{fig:strfactor_ddft} and \ref{fig:strfactor_ddft2} show results for the computed structure factor based on the direct FE results and on the results after the linear mapping to remove artificial correlations, for three different ratios $\Delta x/\ell_0=L/(N_n-1)/\ell_0$  and for both $p1$-elements (Figure \ref{fig:strfactor_ddft}) and $p2$-elements (Figure \ref{fig:strfactor_ddft2}). The numerical results are close to the theoretical curve given by Eq. \eqref{eq:theor_struc_eq2} for all 
discretizations when the linear mapping provided by Equation \eqref{eq:mappingQ} is applied. The direct FE results without the mapping step tend to converge to the theoretical curve too as the mesh is refined for both $p1$- and $p2$-elements, with the results obtained with quadratic elements approximating better the theoretical curve for the structure factor than those obtained with linear elements. However, if no mapping is applied, a finer discretization is needed, as the characteristic length of the mesh needs to be significantly smaller than the typical correlation length to avoid the presence of significant artificial correlations. 
Therefore, applying the decorrelation matrix seems to be advantageous even in the case where there are spatial correlations over a finite length $\ell_0>0$, and its effect is larger for increasing $\Delta x/\ell_0$.

%%%

\section{Discussion and conclusions}\label{sec:conclusions}

In this work, we have presented a numerical approach based on the finite element method to solve stochastic diffusion equations with thermal fluctuations. In particular, we have proposed a general formulation for the discrete fluctuating forcing term that is not restricted to any particular kind of element or shape function, and we have also derived a mass-conserving linear mapping to remove from the discrete solution the artificial correlations introduced by the spatial discretization. 

The proposed mapping provides a solution that has a straightforward physical interpretation regardless of the spatial discretization, so it enables the use of finite element schemes without restrictions. Our analysis makes explicit the fact that the spatial discretization and the coarse-graining are two separate things, and that the artificial correlations due to the spatial discretization depend exclusively on the basis functions used to approximate the finite element solution.
The mapping can be applied to any spatial discretization, including those defined on unstructured meshes.
While the performance of the mapping has been demonstrated for linear equations, the approach is valid for nonlinear equations as well. In fact, the existence of an equivalent mapped solution that is free of artificial correlations can make treating non-linear terms in fluctuating-hydrodynamics equations easier. 

In the case of problems for which the fluctuations are spatially correlated over a finite length, if the spatial discretization is refined enough with respect to the correlation length, the finite element method will eventually provide a solution with physically-meaningful spatial correlations. Nevertheless, the mapping succeeds in removing the artificial correlations related to underresolved fluctuations, thus providing physically meaningful results for coarser meshes. 

\section*{Acknowledgement}
This work was supported by the Ramon y Cajal fellowship RYC2021-030948-I funded by the MCIN/AEI/10.13039/501100011033 and by the EU under the NextGenerationEU/PRTR program, and by the research grant PID2020-113033GB-I00 funded by the MCIN/AEI/10.13039/501100011033.

\bibliography{biblioFH_abbrev.bib} 
\bibliographystyle{ieeetr}

%%%
\appendix

\section{Implementation of the linearized stochastic forcing term}
\label{sec:appendixfluc}

Based on Equation \eqref{eq:cov_b}, an implementation of the stochastic forcing term may be proposed through a decomposition of the diffusion matrix, which for our standard Galerkin discretization is symmetric, and can therefore always be decomposed into
\begin{equation}\label{eq:decompositionD}
\mathbf{D} = \mathbf{A}_D\, \mathbf{A}_D^T \,,
\end{equation}
where matrix $\mathbf{A}_D^T$ is the transpose of matrix $\mathbf{A}_D$.
After integrating in time (see also Section \ref{sec:timeintegrators}), we can compute the array $\mathbf{f}^n$  at every time step $n$ as
\begin{equation}\label{eq:I_linear}
	{f}^n_i = \sqrt{2 u_0}\, a^n_i \,,
\end{equation}
with array $\mathbf{a}^n=\left\{a^n_i\right\}_{i=1}^{N_{\text{dof}}}$ defined as
\begin{equation}\label{eq:I_linear2}
	\mathbf{a}^n = \mathbf{A}_D\,\mathbf{z}^n \,,
\end{equation}
where $\mathbf{z}(t)$ is an array of random numbers of length $N_{\text{dof}}$ following a Gaussian distribution of expected value $0$ and variance $1$.
Matrix $\mathbf{A}_D$ can be computed once before the time-stepping starts, and then $\mathbf{f}^n$ can be computed at each time step based on it. The decomposition given by Equation \eqref{eq:decompositionD} is not unique, but a convenient definition is
\begin{equation}
\mathbf{A}_D=\mathbf{U}_D \sqrt{\boldsymbol{\Sigma}_D} \mathbf{U}_D^T \,,
\label{eq:AD_def}
\end{equation}
where $\mathbf{U}_D$ is a unitary matrix and ${\boldsymbol{\Sigma}_D}$ is a diagonal matrix, satisfying
\begin{equation}
\mathbf{D}=\mathbf{U}_D {\boldsymbol{\Sigma}_D} \mathbf{U}_D^T \,.
\end{equation}
The definition of $\mathbf{A}_D$ provided by Equation \eqref{eq:AD_def} makes $\mathbf{A}_D$ symmetric. Moreover, while in general a matrix $\mathbf{A}_D$ satisfying Equation \eqref{eq:decompositionD} will be a dense matrix, the definition provided by Equation \eqref{eq:AD_def} makes it possible to approximate it by a sparse matrix, which is important to preserve the computational efficiency of the finite element method.
%%%

\section{Dynamic structure factor}
\label{sec:dynamicSF}
In this appendix, we present complementary results to those presented in Sections \ref{sec:results_diff2} and \ref{sec:results_diff2_decorr}. In particular, we present results for the dynamic structure factor, which, unlike the static structure factor discussed in previous sections, contains information about the time evolution of the solution.
 The dynamic structure factor is defined as the Fourier transform of the time-dependent correlation function.
 For the problem discussed in Section \ref{sec:results_diff2}, the analytical expression for the dynamic structure factor is given by
 \begin{equation}
 S_{\text{dyn}}(k,\tau)=\left<\mathcal{U}(k,t) \mathcal{U}(-k,t-\tau)\right>=u_0 e^{-D k^2 \tau}\,.
 \end{equation}

\begin{figure}
\centering
\includegraphics[width=6.7cm]{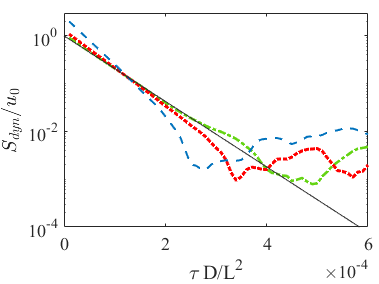}
\includegraphics[width=6.7cm]{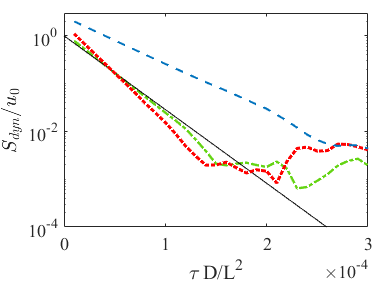}
\caption{Dynamic structure factor of the finite-element solution (before applying the decorrelation matrix) as a function of the time lag for $k=20\,(2\pi/L)$ (left) and $k=30\,(2\pi/L)$ (right): computed with $N_n=51$  (dashed), $N_n=101$  (dotted) and $N_n=201$  (dash-dot), and theoretical dynamic structure factor (solid line). All computations with $p1$-elements, $\Delta t=\num{1e-5}$, $N_t=\num{1e6}$, $u_{0}=10000$, $\alpha=1/2$.}
\label{fig:dynstrfactor_FE}
\end{figure}

 Figure \ref{fig:dynstrfactor_FE} shows the dynamic structure factor corresponding to the finite element solution (before mapping) for two values of the wavenumber $k$ as a function of the time lag $\tau$. As also explained in Section \ref{sec:results_diff2}, for a given value of $k$, the results converge to their theoretical value as the mesh is refined. 

\begin{figure}
\centering
\includegraphics[width=6.7cm]{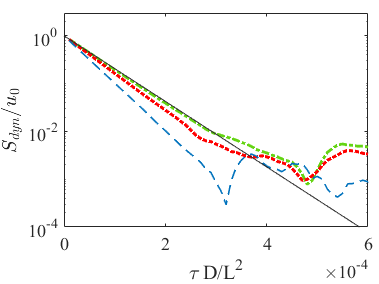}
\includegraphics[width=6.7cm]{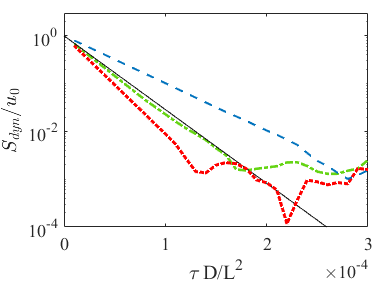}
\caption{Dynamic structure factor of the mapped solution (after applying the decorrelation matrix) as a function of the time lag for $k=20\,(2\pi/L)$ (left) and $k=30\,(2\pi/L)$ (right): computed with $N_n=51$  (dashed), $N_n=101$  (dotted) and $N_n=201$  (dash-dot), and theoretical dynamic structure factor (solid line). All computations with $p1$-elements, $\Delta t=\num{1e-5}$, $N_t=\num{1e6}$, $u_{0}=10000$, $\alpha=1/2$.}
\label{fig:dynstrfactor_mapped}
\end{figure}

 Figure \ref{fig:dynstrfactor_mapped} displays the dynamic structure factor corresponding to the mapped finite element solution for the same two values of $k$. A key difference with respect to the finite element results before mapping (Figure \ref{fig:dynstrfactor_FE}) is that the mapped results present a good agreement with the theoretical value for $\tau$ values close to zero regardless of the size of the mesh. For larger values of $\tau$, the dynamic structure factor of the mapped results converges as well towards the theoretical one as the mesh is refined.
 
\begin{figure}
\centering
\includegraphics[width=6.7cm]{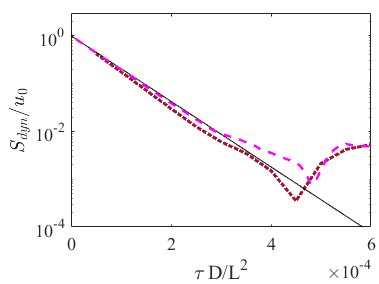}
\caption{Influence of the time step on the dynamic structure factor of the mapped solution for $k=20\,(2\pi/L)$: computed with $\Delta t=\num{1e-5}$  (dashed) and $\Delta t=\num{5e-5}$ (dotted), and theoretical dynamic structure factor (solid line). All computations with $p1$-elements, $N_n=201$, $N_t \Delta t=10 L^2/D$, $u_{0}=10000$, $\alpha=1/2$.}
\label{fig:dynstrfactor_mapped_dt}
\end{figure}

\begin{figure}
\centering
\includegraphics[width=6.7cm]{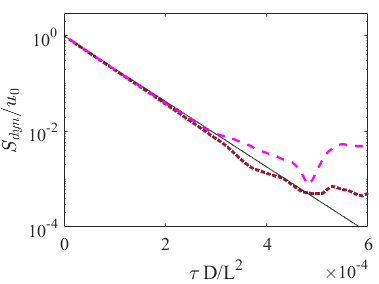}
\caption{Influence of the total simulation time on the dynamic structure factor of the mapped solution for $k=20\,(2\pi/L)$: computed with $N_t \Delta t =10 L^2/D$  (dashed) and $N_t \Delta t=100 L^2/D$ (dotted), and theoretical dynamic structure factor (solid line). All computations with $p1$-elements, $N_n=201$, $\Delta t=\num{1e-5}$, $u_{0}=10000$, $\alpha=1/2$.}
\label{fig:dynstrfactor_mapped_Nt}
\end{figure}

Figures \ref{fig:dynstrfactor_mapped_dt} and \ref{fig:dynstrfactor_mapped_Nt} show the effect of the time step $\Delta t$ and of the total simulated time $N_t \Delta t$ on the computed structure factor based on the mapped solution, and illustrate the convergence of the numerical solution towards the theoretical one as the time step decreases and the total simulated time increases. For values of $\tau$ below a threshold, decreasing the time step leads to a better agreement with the theoretical solution, as seen in Figure \ref{fig:dynstrfactor_mapped_dt}. For values of $\tau$ above a certain threshold, the value of the dynamic structure factor is small, and the residual error due to the finite time $N_t \Delta t$ that is simulated dominates the solution. As shown in Figure  \ref{fig:dynstrfactor_mapped_Nt}, this error decreases as the total simulated time $N_t \Delta t$ increases. 

These results show that the numerical approach based on applying the decorrelation matrix to the finite element results, as proposed in Section \ref{sec:decorrelation}, can capture the time evolution of the solution.

\end{document}